\documentclass[a4paper,12pt]{amsart}
\usepackage{amsfonts}
\usepackage{amssymb}
\usepackage{ifthen}
\usepackage{graphicx}
\nonstopmode \numberwithin{equation}{section}
\usepackage[usenames]{color}
\newtheorem{thm}{Theorem}[section]
\newtheorem{lem}{Lemma}[section]
\newtheorem{cor}{Corollary}[section]

\newtheorem{cl}{Claim}[section]
\newtheorem{ca}{Case}[section]
\newtheorem{sca}{Subcase}[section]
\newtheorem{scl}[section]{Subclaim}
\newtheorem{conj}[equation]{Conjecture}

\theoremstyle{definition}
\newtheorem{defn}{Definition}[section]
\newtheorem{op}[equation]{Open Problem}
\newtheorem{ques}[equation]{Question}
\newtheorem{rem}{Remark}[section]
\newtheorem{exam}[equation]{Example}

\newcounter {own}
\def\theown {\thesection       .\arabic{own}}

\newenvironment{pf}[1][]{%
 \vskip 3mm
 \noindent
 \ifthenelse{\equal{#1}{}}%
  {{\slshape Proof. }}%
  {{\slshape #1.} }%
 }%
{\qed\bigskip}

\newcounter{alphabet}
\newcounter{tmp}
\newenvironment{Thm}[1][]{\refstepcounter{alphabet}%
\bigskip%
\noindent%
{\bf Theorem \Alph{alphabet}}%
\ifthenelse{\equal{#1}{}}{}{ (#1)}%
{\bf .} \itshape}{\vskip 8pt}

\makeatletter
\renewcommand{\Ref}[1]{\@ifundefined{r@#1}{}{\setcounter{tmp}{\ref{#1}}\Alph{tmp}}}
\makeatother

\newenvironment{Lem}[1][]{\refstepcounter{alphabet}%
\bigskip%
\noindent%
{\bf Lemma \Alph{alphabet}}%
{\bf .} \itshape}{\vskip 8pt}

\newcommand{\IR}{{\mathbb R}}

\newcommand{\IC}{{\mathbb C}}

\newcommand{\diam}{{\operatorname{diam}}}

\def\be{\begin{equation}}
\def\ee{\end{equation}}

\newcommand{\bee}{\begin{enumerate}}
\newcommand{\eee}{\end{enumerate}}

\newcommand{\blem}{\begin{lem}}
\newcommand{\elem}{\end{lem}}
\newcommand{\bthm}{\begin{thm}}
\newcommand{\ethm}{\end{thm}}
\newcommand{\bcor}{\begin{cor}}
\newcommand{\ecor}{\end{cor}}
\newcommand{\beg}{\begin{exam}}
\newcommand{\eeg}{\end{exam}}
\newcommand{\begs}{\begin{examples}}
\newcommand{\eegs}{\end{examples}}
\newcommand{\bdefe}{\begin{defn}}
\newcommand{\edefe}{\end{defn}}
\newcommand{\bprob}{\begin{prob}}
\newcommand{\eprob}{\end{prob}}
\newcommand{\bques}{\begin{ques}}
\newcommand{\eques}{\end{ques}}
\newcommand{\bei}{\begin{itemize}}
\newcommand{\eei}{\end{itemize}}
\newcommand{\bcon}{\begin{conj}}
\newcommand{\econ}{\end{conj}}
\newcommand{\bop}{\begin{op}}
\newcommand{\eop}{\end{op}}

\newcommand{\bca}{\begin{ca}}
\newcommand{\eca}{\end{ca}}
\newcommand{\bsca}{\begin{sca}}
\newcommand{\esca}{\end{sca}}

\newcommand{\bcl}{\begin{cl}}
\newcommand{\ecl}{\end{cl}}

\newcommand{\bscl}{\begin{scl}}
\newcommand{\escl}{\end{scl}}

\newcommand{\bcons}{\begin{conjs}}
\newcommand{\econs}{\end{conjs}}
\newcommand{\bprop}{\begin{propo}}
\newcommand{\eprop}{\end{propo}}
\newcommand{\br}{\begin{rem}}
\newcommand{\er}{\end{rem}}
\newcommand{\brs}{\begin{rems}}
\newcommand{\ers}{\end{rems}}
\newcommand{\bo}{\begin{obser}}
\newcommand{\eo}{\end{obser}}
\newcommand{\bos}{\begin{obsers}}
\newcommand{\eos}{\end{obsers}}
\newcommand{\bpf}{\begin{pf}}
\newcommand{\epf}{\end{pf}}
\newcommand{\ba}{\begin{array}}
\newcommand{\ea}{\end{array}}
\newcommand{\beq}{\begin{eqnarray}}
\newcommand{\beqq}{\begin{eqnarray*}}
\newcommand{\eeq}{\end{eqnarray}}
\newcommand{\eeqq}{\end{eqnarray*}}

\newcounter{minutes}
\divide\time by 60
\newcounter{hours}
\multiply\time by 60 \addtocounter{minutes}{-\time}

\begin{document}

\bibliographystyle{amsplain}
\title{A note on V\"ais\"al\"a's problem concerning free quasiconformal mappings}
\thanks{$^\dagger$ File:~\jobname .tex,
          printed: \number\year-\number\month-\number\day,
          \thehours.\ifnum\theminutes<10{0}\fi\theminutes}

\author{Qingshan Zhou}
\address{Qingshan Zhou, School of Mathematics and Big Data, Foshan University,  Foshan, Guangdong 528000, People's Republic of China} \email{qszhou1989@163.com; q476308142@qq.com}

\author{Yaxiang  Li}
\address{Yaxiang Li,  Department of Mathematics, Hunan First Normal University, Changsha,
Hunan 410205, People's Republic of China} \email{yaxiangli@163.com}

\author{Antti Rasila${}^{~\mathbf{*}}$}
\address{Antti Rasila,  Technion -- Israel Institute of Technology, Guangdong Technion, 241 Daxue Road, Shantou, Guangdong 515063, People's Republic
of China} \email{antti.rasila@gtiit.edu.cn; antti.rasila@iki.fi}

\date{}
\keywords{Quasihyperbolic metric, free quasiconformal mapping, quasisymmetric mapping.\\
${}^{\mathbf{*}}$ Corresponding author}

\begin{abstract} In this paper, we provide partial solutions to a problem raised by V\"ais\"al\"a on local properties of free quasiconformal mappings. In particular, we show that a locally free quasiconformal mapping is globally free quasiconformal under the condition of locally relative quasisymmetry.

\end{abstract}

\thanks{}

\maketitle\pagestyle{myheadings} \markboth{}{A note on V\"ais\"al\"a's problem}

\section{Introduction and main results}\label{sec-1}
Many results of conformal mapping theory have their counterparts in the setting of quasiconformal mappings in $n$-dimensional Euclidean spaces. For instance, Gehring and his students Palka and Osgood \cite{GO,GP} established  a Schwarz-Pick type theorem for quasiconformal mappings in ${\mathbb R}^n$ $(n\ge 2)$ by introducing the quasihyperbolic metric (see Definition \ref{i-2}). This metric is a generalization of the classical hyperbolic metric in the unit disk. Since its first appearance, the quasihyperbolic metric has become an important tool in the geometric function theory of Euclidean spaces, and it is the basic element of the definition of free quasiconformality.

Based on results of Gehring and Osgood, Tukia and V\"ais\"al\"a \cite{TV81} considered another class of homeomorphisms, solid mappings, which forms an intermediate class between quasiconformal and quasisymmetric mappings. Solid mappings are very useful to the study of the equivalence of these two classes of mappings in  different settings (see e.g. \cite{HL, HRWZ,Vai-2}).

From late 1980's onwards, V\"ais\"al\"a has developed the theory of (dimension) free quasiconformal mappings in Banach spaces, as an attempt to further generalize the higher dimensional quasiconformal mapping theory, see \cite{Vai-1, Vai-2, Vai-3, Vai-4, Vai-5}. The main advantage of this approach is that it avoids using volume integrals and conformal modulus, allowing one to study the quasiconformality of mappings in infinite dimensional Banach spaces and metric spaces without doubling volume measure. This line of research has recently attracted substantial interest, and some recent results can be found from \cite{HLL,HPRW,LPV,rt2} and references therein.

In this paper, we focus on local properties of free quasiconformal mappings in Banach spaces. In particular, we investigate several local conditions under which the free quasiconformality of a mapping can be established.

Throughout this paper, following notations and terminology of \cite{Vai-5}, we assume that $E$ and $E'$ are real Banach
spaces with dimension at least $2$, and that $D\subsetneq E$ and $D'\subsetneq E'$ are domains (open and connected sets). The norm of a vector $z$ in $E$ is written as $|z|$, and for each pair of points $z_1$, $z_2$ in $E$, the distance between them is denoted by $|z_1-z_2|$, and the closed line segment with endpoints $z_1$ and $z_2$ is $[z_1, z_2]$. The distance from $z\in D$ to the boundary $\partial D$ of $D$ is denoted by $d_D(z)$. For a bounded set $A$ in $E$, $\diam (A)$ means the diameter of $A$. Open and closed balls with center $x$ and radius $r$ are denoted by ${B}(x,r)$ and $\overline{B}(x,r)$, respectively. The boundary of ${B}(x,r)$ is
denoted by ${S}(x,r)$.

\bdefe \label{def1.7} Let $\varphi:[0,\infty)\to [0,\infty)$ be a homeomorphism with $\varphi(t)\geq t$. We say
that a homeomorphism $f: D\to D'$ is {\it $\varphi$-semisolid} if
$$k_{D'}(f(x),f(y))\leq \varphi(k_D(x,y))$$
for all $x$, $y\in D$, and {\it $\varphi$-solid} if both $f$ and $f^{-1}$
satisfy this condition. We say that $f$ has fully a given property if for each subdomain $G\subset D$, $f|_G$ has this property. For example, a homeomorphism $f: D\to D'$ is called {\it fully $\varphi$-semisolid} if for each subdomain $G\subset D$, $f|_G$ is $\varphi$-semisolid. Note that the classes of fully $\varphi$-semisolid maps, fully $\varphi$-solid maps,
 fully $(M,C)$-CQH maps, and fully $(M, C)$-CDR maps are defined in a similar way. For the definitions of $(M,C)$-CQH maps and $(M, C)$-CDR maps see Definition \ref{i-1}.
 \edefe

 Fully $\varphi$-solid mappings are also called free
$\varphi$-quasiconformal mappings, or $\varphi$-FQC mappings. In \cite{TV81}, Tukia and V\"ais\"al\"a demonstrated that a homeomorphism between two proper domains of $\mathbb{R}^n$ is $K$-quasiconformal if and only if it is $\varphi$-FQC, where $K$ and $\varphi$ depend only on each other and $n$.

Note that the quasiconformality is equivalent to local quasiconformality, that is, the mapping is quasiconformal on a neighborhood of every point $x$ in the domain $D$. It is natural to ask whether this result holds for free quasiconformal mappings in the setting of Banach spaces. Indeed, V\"ais\"al\"a \cite{Vai-1} has raised the following question:

\bop\label{Con2}$($\cite[Section 7]{Vai-1}$)$ Suppose that $f: D \to D'$ is a
homeomorphism and that each point has a neighborhood $G \subseteq D$ such that $f_G: G \to f(G)$ is $\varphi$-FQC.
Is $f$ $\varphi'$-FQC with $\varphi'=\varphi'(\varphi)$? \eop

%

Denote $\overline{B}_x=\overline{B}(x,d_D(x))$ and $S_x=S(x,d_D(x))$. We shall provide a partial solution to Problem \ref{Con2} under the additional assumption that $f|_{\overline{B}_x}$ is quasisymmetric relative  to $S_x$ (briefly QS relative to $S_x$) for all $x\in D$.  For the definition of relative quasisymmetric mappings see Definition \ref{def1-0}. Our main result is the following:

\begin{thm}\label{thm2.0} Let $\varphi,\eta:[0,\infty)\to [0,\infty)$ be homeomorphisms. Suppose that $f: D \to D'$ is a homeomorphism. If for each point $x\in D$,
\begin{enumerate}
  \item $f|_{B_x}$ is $\varphi$-FQC, and
  \item $f|_{\overline{B}_x}$ is $\eta$-QS relative to $S_x$,
\end{enumerate}
then $f$ is $\psi$-FQC with $\psi=\psi(\varphi,\eta)$.
\end{thm}

\br If $f|_{B_x}$ is $\eta$-QS for all $x\in D$, then it satisfies the conditions $(1)$ and $(2)$ of Theorem \ref{thm2.0}. This can be obtained by using \cite[Theorem 7.9]{Vai-5} and \cite[Theorem 6.12]{Vai-5}. That is to say, a mapping that is quasisymmetric on all Whitney type balls $B_x$, is a globally free quasiconformal mapping. In a sense, this condition is also necessary. Indeed, it follows from \cite[Theorem 7.9]{Vai-5} that a globally FQC mapping is quasisymmetric on every Whitney type ball strictly contained in the domain.
\er

Then we give an example to show that free quasiconformality does not imply the relative quasisymmetry. Recall that the concept of uniform domains was introduced by
 Martio and Sarvas \cite{MS} in connection to global injectivity properties of locally injective mappings.  This class of domains has numerous geometric and function theoretic properties that are useful in many fields of mathematical analysis (see e.g. \cite{GP,HPRW,Vai-2,Vai-5,ZRL}).

\beg
Let $D=\{z=x+iy\in \IC:\; |z|<1\}$, $D'=D\setminus\{z=x+iy\in \IC:\; y=0,0<x<1 \}$, and let $f:$ $D\to D'=f(D)$ be a conformal mapping. Then $f$ is $\varphi$-FQC, but not relative quasisymmetric. Otherwise, by \cite[Theorem 1.1]{HLVW} it would follow that $D'$ is a uniform domain (see Definition \ref{def1.1}), which is a contradiction.
\eeg

By using Theorem \ref{thm2.0}, we further investigate other (local) conditions such that the locally free quasiconformality implies globally free quasiconformality. The following result can be viewed as a second partial solution to Problem \ref{Con2}.

\begin{thm}\label{cor1} Let $\varphi:[0,\infty)\to [0,\infty)$ be a homeomorphism, and let $L, c\geq 1$. Suppose that $f: D \to D'$ is a homeomorphism. If for each point $x\in D$,
\begin{enumerate}
  \item $f|_{B_x}$ is $\varphi$-FQC, and
  \item $\diam (f(B_x))\leq L d_{f(B_x)}(f(x))$ and $f(B_x)$ is $c$-uniform,
\end{enumerate}
then $f$ is $\psi$-FQC  with $\psi=\psi(\varphi,L,c)$.
\end{thm}

Next, we study the local properties of FQC mappings satisfying the coarse distance ratio (CDR) distortion property on all balls $B_x$ of the domain. Note that the CDR mappings (see Definition \ref{i-1} (\ref{i-1b})) were introduced by the authors in \cite{ZRL}, where we also studied the relationships between the classes of CDR and FQC mappings and established invariance results for certain classes of domains under CDR mappings.
As an application of Theorem \ref{cor1}, we obtain the following:

\begin{cor}\label{cor2}  Let $\varphi:[0,\infty)\to [0,\infty)$ be a homeomorphism, and let $L,M\geq 1$, $C\geq 0$. Suppose that $f: D \to D'$ is a homeomorphism. If for each point $x\in D$,
\begin{enumerate}
  \item $f|_{B_x}$ is $\varphi$-FQC,
  \item $\diam (f(B_x))\leq L d_{f(B_x)}(f(x))$ and $f|_{B_x}$ is $(M,C)$-CDR,
\end{enumerate}
then $f$ is $\psi$-FQC  with $\psi=\psi(\varphi,L,M,C)$.
\end{cor}

Generally, a FQC mapping is not CDR.

\beg $($\cite[Example 3.4]{ZRL}$)$
Let $D=\{z=x+iy\in \IC:\; |z|<1\}$, $D'=\{z=x+iy\in \IC:\; 0<y<1 \}$, and let $f:$ $D\to D'=f(D)$ be a conformal mapping. Then $f$ is $\varphi$-FQC, but not CDR.
\eeg

The rest of this paper is organized as follows.  In Section \ref{sec-2}, we recall necessary definitions and preliminary results.  Section \ref{sec-3} is devoted to the proofs of Theorems \ref{thm2.0}, \ref{cor1} and Corollary \ref{cor2}.

\section{Preliminaries}\label{sec-2}

\subsection{Quasihyperbolic metric and distance ratio metric} In this part, we assume that $E$ and $E'$ are real Banach spaces with dimension at least $2$, and that $D\subsetneq E$ and $D'\subsetneq E'$ are domains. The quasihyperbolic metric of domains was introduced by Gehring with his students Palka and Osgood \cite{GO,GP}. Recall the following definition:

\bdefe\label{i-2} Let $x,y\in D$. The {\it quasihyperbolic metric} $k_D$ of $D$ is defined by
$$k_D(x,y)=\inf \int_\gamma \frac{|dz|}{d_D(z)},$$
where the infimum is taken over all rectifiable curves $\gamma$ joining $x$ and $y$ in $D$.
\edefe

The quasihyperbolic metric is closely connected to the distance ratio metric, defined as follows:

\bdefe\label{i-3} Let $x,y\in D$. The {\it distance ratio metric} $j_D$ of $D$ is defined by
$$j_D(x,y)=\log\Big(1+\frac{|x-y|}{\min\{d_{D}(x), d_{D}(y)\}}\Big) .$$
\edefe

Note that the distance ratio metric always gives the elementary lower bound $j_D(x,y) \le k_D(x,y)$ for the quasihyperbolic distance. Next we recall the definition of uniform domains, which were originally introduced by Martio and Sarvas \cite{MS}:

\bdefe\label{def1.1} Let $c\geq 1$. We say that $D$ is $c$-{\it uniform} if each pair of points $x$, $y$ in $D$ can be connected by an arc $\gamma$ in $D$ satisfying:
\begin{enumerate}
\item $\ell(\gamma)\leq c\,|x-y|$, and
\item $\min\{\ell(\gamma[x,z]),\ell(\gamma[z,y])\}\leq c\,d_D(z)$ for all $z\in \gamma$,
\end{enumerate}
\noindent where $\ell(\gamma)$ denotes the length of $\gamma$, $\gamma[x,z]$ the part of $\gamma$ between $x$ and $z$.
\edefe

In particular, it is well known that balls are uniform domains in Banach spaces:

\begin{Lem}\label{lem-v} $($\cite[Theorem 6.5]{Vai-2}$)$
Every ball in a Banach space is $c_0$-uniform with a universal constant $c_0$.
\end{Lem}

Moreover, we need the following fact on uniform domains.

\begin{Thm}\label{lem-z} {\rm (\cite[Theorem 6.16]{Vai-3})} Let $D\subsetneq E$ be a domain. Then $D$ is $c$-uniform if and only if $k_D(x,y)\leq c_1j_D(x,y)+c_2$ for all $x,y\in D$.
\end{Thm}

In the above, the additive constant $c_2$ can be absorbed into the constant $c_1$. See, for example, \cite[Theorem 6.16]{Vai-2}. Finally, we introduce the definitions of two classes of coarse maps.

\bdefe\label{i-1}
Let $M\geq 1$ and $C\geq 0$, and let $f:\;D\to D'$ be a homeomorphism.
\begin{enumerate}
\item\label{i-1a} $f$ is called {\it $C$-coarsely $M$-quasihyperbolic} ($(M,C)$-CQH) if for all $x,y\in D$,
$$M^{-1} k_{D}(x,y)-C \leq k_{D'}(f(x),f(y))\leq  {M}k_{D}(x,y) +C.$$

\item\label{i-1b} $f$ is  {\it $C$-coarsely $M$-distance ratio} ($(M,C)$-{\it CDR}) if for all $x,y\in D$,
$$M^{-1} j_{D}(x,y)-C \leq j_{D'}(f(x),f(y))\leq  {M}j_{D}(x,y) +C.$$

\end{enumerate}
\edefe

\subsection{Quasisymmetric mappings and quasim\"obius mappings}

Let $(X,|\cdot|)$ be a metric space.
By a triple in $X$ we mean an ordered sequence $T=(x,y,z)$ of three
distinct points in $X$. The ratio of $T$ is the number
$\rho(T)=|y-x|/|z-x|.$ If $f: X\to Y$ is  an injective
map, the image of a triple  $T=(x,y,z)$  is the triple
$f(T)=(f(x),f(y),f(z))$.

Suppose that  $A\subseteq X$. A triple  $T=(x,y,z)$ in $X$ is said to
be a triple in the pair $(X, A)$ if $x\in A$ or if $\{y,z\}\subseteq
A$. Equivalently, both $|y-x|$ and $|z-x|$ are distances from a
point in $A$.

\bdefe \label{def1-0} Let $X$ and $Y$ be two metric spaces, and let
$\eta: [0, \infty)\to [0, \infty)$ be a homeomorphism. Suppose
$A\subseteq X$. An embedding $f: X\to Y$ is said to be {\it
$\eta$-quasisymmetric} relative to $A$ ($\eta$-$QS$ relative to
$A$) if $\rho(f(T))\leq \eta(\rho(T))$  for each triple $T$ in
$(X,A)$. Moreover, $\eta$-quasisymmetric relative to $X$ is equivalent to usual $\eta$-quasisymmetric ($\eta$-QS).\edefe

\bdefe \label{def1-0'} Let  $0<q<1$. Assume that $D\subsetneq E$ and $D'\subsetneq E'$ are domains in real Banach spaces $E$ and $E'$, respectively.  A homeomorphism $f:D\to D'$ is
{\it $q$-locally $\eta$-quasisymmetric} ($q$-locally $\eta$-QS) if $f|_{{B}(x,qr)}$ is
$\eta$-QS for all $x\in D$ and $0<r\leq d_{D}(x)$. \edefe



 It is known that each $K$-QC mapping in $\IR^n$ is $q$-locally $\eta$-QS for every $0<q<1$ with $\eta=\eta(K, q, n)$. Conversely, each $q$-locally $\eta$-QS mapping in $\IR^n$ is a $K$-QC mapping with $K=\eta(1)^{n-1}$ by the metric definition of quasiconformality (cf. \cite[5.6]{Vai-1}). In \cite{Vai-5}, V\"ais\"al\"a further proved the following:

\begin{Thm}\label{Thm2'}$($\cite[Theorem 7.9]{Vai-5}$)$ Assume that $D\subsetneq E$ and $D'\subsetneq E'$ are domains in real Banach spaces $E$ and $E'$, respectively.
For a homeomorphism $f: D\to D'$, the following conditions are quantitatively equivalent:

\bee
\item $f$ is $\varphi$-FQC;
\item $f$ and $f^{-1}$ are both fully $(M,C)$-CQH;
\item $f$ and $f^{-1}$ are  both $q$-locally $\eta$-QS;
\item For every $0<q<1$, there is some  $\eta_q$ such that both $f$ and $f^{-1}$ are  $q$-locally $\eta_q$-QS;
\eee\end{Thm}

Let $(X,|\cdot|)$ be a metric space. A quadruple in $X$ is an ordered sequence $Q=(x,y,z,w)$ of four
distinct points in $X$. The cross ratio of $Q$ is defined to be the
number
$$\tau(Q)=\tau(x,y,z,w)=\frac{|x-z|}{|x-y|}\cdot\frac{|y-w|}{|z-w|}.$$

Denote $\dot{X}=X\cup \{\infty\}$. Observe that the definition of cross ratio can be extended in
the well known manner to the case where one of the points is
$\infty$. For example,
$$\tau(x,y,z,\infty)= \frac{|x-z|}{|x-y|}.$$ If $X_0 \subseteq \dot{X}$ and if $f: X_0\to \dot{Y}$
is an injective map, the image of a quadruple $Q$ in $X_0$ is the
quadruple $f(Q)=(f(x),f(y),f(z),f(w))$. Suppose that $A\subseteq X_0$. We say
that a quadruple $Q=(x,y,z,w)$ in $X_0$ is a quadruple in the pair
$(X_0, A)$ if $\{x,w\}\subseteq A$ or $\{y,z\}\subseteq A$.
Equivalently, all four distances in the definition of $\tau(Q)$ are
(at least formally) distances from a point in $A$.

\bdefe \label{def2'} Let ${X}$ and ${Y}$ be two metric
spaces and let $\eta: [0, \infty)\to [0, \infty)$ be a
homeomorphism. Suppose $A\subseteq \dot{X}=X\cup\{\infty\}$. An embedding $f:
\dot{X}\to \dot{Y}$ is said to be {\it $\eta$-quasim\"obius
relative to $A$}, or $\eta$-$QM$ relative to $A$, if the inequality
$\tau(f(Q))\leq \eta(\tau(Q))$ holds for each quadruple in $(X,A)$.
\edefe
Apparently, $\eta$-$QM$ relative to $X$ is equivalent to ordinary $\eta$-quasim\"obius ($\eta$-QM).

\section{Proofs of Theorem \ref{thm2.0}, Theorem \ref{cor1}, and Corollary \ref{cor2}}\label{sec-3}

\subsection{Proof of Theorem \ref{thm2.0}}Fix $x_0\in D$. Note that $B_x=B(x,d_D(x))$ for all $x\in D$. We first claim that
\be\label{1.0} d_{D'}(f(x_0))\leq 2\eta(1)d_{f(B_{x_0})}(f(x_0)).\ee

Let $\varepsilon>0$ be sufficiently small so that $\eta(1)\eta(\varepsilon)\leq \frac{1}{2}$, and let $\varepsilon_i=2^{-i} \varepsilon$ for non-negative integers $i$. First, we see from the arguments in the proof of \cite[Theorem 3.10]{Vai-2} that there is a  sequence of points $\{x_i\}_{i\in \mathbb{N}}$ and an end-cut $\gamma$ of $D$ from $x_0$ such that the following hold:
\begin{enumerate}
\item\label{1.1}
$\gamma_n=\bigcup_{i=0}^{n-1}[x_i,x_{i+1})$ and $\gamma=\bigcup_{n=0}^{\infty}\gamma_n$, where $[x_i,x_{i+1}]$ is a line segment in $\overline{D}$;
\item\label{1.2}
$x_{i+1}\in S(x_i,d_D(x_i))$ and $d_D(x_{i+1})\leq \varepsilon_{i+1} d_D(x_i)$;
\item\label{1.3}
$\lim\limits_{n\to\infty}x_n=a\in\partial D$.
\end{enumerate}

Without loss of generality, we may assume that $x_i\in D$ for all $i\in \mathbb{N}$. We  show that
\be\label{1.4} |f(x_{i+1})-f(x_{i+2})|\leq \frac{1}{2}|f(x_{i})-f(x_{i+1})|\ee
for all $i\in \mathbb{N}$. This can be seen as follows. Fix $i\in \mathbb{N}$. Choose a point $z_{i+1}\in [x_i,x_{i+1}]\cap {S}(x_{i+1},d_D(x_{i+1}))$. We check that $z_{i+1}\neq x_{i+2}$. Indeed, we see from the property (\ref{1.2}) of the sequence $\{x_n\}$ that
\beq \nonumber d_D(x_{i+2})&\leq& \varepsilon_{i+2} d_D(x_{i+1})= \varepsilon_{i+2} |x_{i+1}-x_{i+2}|
\\ \nonumber &=&  \varepsilon_{i+2} |x_{i+1}-z_{i+1}|= \varepsilon_{i+2}d_{B_{x_i}}(z_{i+1})
\\ \nonumber &\leq& \varepsilon_{i+2}d_{D}(z_{i+1}),
\eeq
as required.

Because $f$ is $\eta$-QS relative to $S(x_{i+1},d_D(x_{i+1}))$ on the closed ball $\overline{B}(x_{i+1},d_D(x_{i+1}))$, we have
$$\frac{|f(x_{i+2})-f(x_{i+1})|}{|f(z_{i+1})-f(x_{i+1})|}\leq \eta\Big(\frac{|x_{i+2}- x_{i+1}| }{|z_{i+1}- x_{i+1}|}\Big)=\eta(1).$$
Similarly, because $f$ is $\eta$-QS relative to $S(x_{i},d_D(x_{i}))$ on the closed ball $\overline{B}(x_{i},d_D(x_{i}))$, we obtain
$$\frac{|f(x_{i+1})-f(z_{i+1})|}{|f(x_{i+1})-f(x_{i})|}\leq \eta\Big(\frac{|x_{i+1}- z_{i+1}| }{|x_{i+1}- x_i|}\Big) \leq \eta(\varepsilon_{i+1}),$$
by the choice of $z_{i+1}$. From this it follows that
\beq \nonumber |f(x_{i+1})-f(x_{i+2})|&\leq& \eta(1)|f(x_{i+1})-f(z_{i+1})|
\\ \nonumber &\leq&  \eta(1)\eta(\varepsilon)|f(x_{i})-f(x_{i+1})|
\\ \nonumber &\leq& \frac{1}{2}|f(x_{i})-f(x_{i+1})|,
\eeq
which yields (\ref{1.4}).

Moreover, by (\ref{1.4}), we obtain
$$|f(x_0)-f(x_{n+1})|\leq \sum_{i=0}^n |f(x_{i})-f(x_{i+1})| \leq 2 |f(x_{0})-f(x_{1})|.$$
Furthermore, from the property (\ref{1.3}) of the sequence $\{x_n\}$ it follows that
$$d_{D'}(f(x_0))\leq |f(x_0)-f(a)|\leq 2|f(x_0)-f(x_1)|\leq 2\eta(1)d_{f(B_{x_0})}(f(x_0)).$$
This proves (\ref{1.0}).

Next, we show that both $f$ and its inverse map $f^{-1}$ are locally $\theta$-QS for some homeomorphism $\theta:[0,\infty)\to [0,\infty)$. Because $f: B_{x} \to f(B_{x})$ is $\varphi$-FQC, it follows from Theorem \Ref{Thm2'} that $f$ is $\theta$-quasisymmetric on $B(x,\frac{1}{2}d_{B_x}(x))$, and similarly, $f^{-1}$ is $\theta$-quasisymmetric on the ball $B(f(x),\frac{1}{2} d_{f(B_x)}(f(x)))$ with $\theta=\theta(\varphi)$. Note that because $d_D(x)= d_{B_x}(x)$, one immediately finds that $f$ is $\theta$-QS on $B(x,\frac{1}{2}d_D(x))$. Then by (\ref{1.0}), we have
$$B\Big(f(x),\frac{1}{4\eta(1)} d_{D'}(f(x))\Big) \subseteq B\Big(f(x),\frac{1}{2} d_{f(B_x)}(f(x))\Big)$$
and thus $f^{-1}$ is $\theta$-QS on $B\big(f(x),\frac{1}{4\eta(1)} d_{D'}(f(x))\big)$, as desired. Hence, Theorem \Ref{Thm2'} ensures that $f$ is $\psi$-FQC with $\psi=\psi(\varphi,\eta)$.
\qed

\subsection{Proof of Theorem \ref{cor1} }
By Theorem \ref{thm2.0}, we only need to show that $f|_{\overline{B}_x}$ is $\eta$-QS relative to $S_x$ for all $x\in D$.

Fix $x\in D$. It follows from Theorem \Ref{Thm2'} that $f|_{{B_x}}$ is $(M,C)$-CQH for some constants $M\geq 1$ and $C\geq 0$ depending only on $\varphi$. Then by Lemma \Ref{lem-v} we see that $B_x$ is $c_0$-uniform with a universal constant $c_0$. Because $f(B_x)$ is $c$-uniform, by \cite[Theorem 11.8]{Vai-5} we know that $f|_{{B_x}}$ extends to a homeomorphism $f|_{\overline{B}_x}$ which is $\theta$-QM relative to $S_x$ with $\theta=\theta(C,M,c,c_0)$.
Because $\diam f(B_x)\leq L d_{f(B_x)}(f(x))$, we observe from \cite[Theorem 6.31]{Vai-5} that $f|_{\overline{B}_x}$ is  $\eta$-QS relative to $S_x$ with $\eta=\eta(\theta,L)$.
\qed

\subsection{Proof of  Corollary  \ref{cor2}} By Theorem \ref{cor1}, it suffices to show that for all $x\in D$, $f(B_x)$ is $c'$-uniform for some constant $c'\geq 1$.

On the one hand, it follows from Lemma \Ref{lem-v} that $B_x$ is $c_0$-uniform with a universal constant $c_0$. Thus, by Theorem \Ref{lem-z}, we know that there are universal constants $c_1$ and $c_2$ such that for all $y,z\in B_x$,
\be\label{z-1} k_{B_x}(y,z)\leq c_1j_{B_x}(y,z)+c_2.\ee

On the other hand, we see from Theorem \Ref{Thm2'}  that $f|_{{B_x}}$ is $(M',C')$-CQH for some constants $M'\geq 1$ and $C'\geq 0$ depending only on $\varphi$. Because $f|_{B_x}$ is $(M,C)$-CDR, for  all $y,z\in B_x$ we thus have by (\ref{z-1}) that
\beq \nonumber k_{f(B_x)}(f(y),f(z))&\leq& M'k_{B_x}(y,z)+C'
\\ \nonumber &\leq&  M'(c_1j_{B_x}(y,z)+c_2)+C'
\\ \nonumber &\leq&  M'(c_1(Mj_{f(B_x)}(f(y),f(z))+C)+c_2)+C'.
\eeq
This, together with Theorem \Ref{lem-z} shows that $f(B_x)$ is $c'$-uniform for some $c'\geq 1$.
\qed

\bigskip

{\bf Funding:} The first author was supported by NNSF of
		China (No. 11901090), and by Department of Education of Guangdong Province, China (No. 2018KQNCX285). The second author was supported by NNSF of	 China (No. 11601529, 11671127) and Hunan Provincial Key Laboratory of Mathematical Modeling and Analysis in Engineering  (No. 2017TP1017). The third author was supported by NNSF of
		China (No. 11971124) and NSF of Guangdong Province.

\end{document}